\documentclass{amsart}
\usepackage{amscd, amsmath, amsthm, amssymb}

\newtheorem{theorem}{Theorem}[section]
\newtheorem{lemma}[theorem]{Lemma}
\newtheorem{proposition}[theorem]{Proposition}

\theoremstyle{definition}     
\newtheorem{definition}[theorem]{Definition}

\newtheorem{question}[theorem]{Question}

\theoremstyle{remark}
\newtheorem{remark}[theorem]{Remark}

\numberwithin{equation}{section}

\begin{document}

\title[hyperk\"aher manifolds ]
{Mordell-Weil groups of a hyperk\"ahler manifold 
- a question of F. Campana}

\author[K. Oguiso]{Keiji Oguiso}
\address{Department of Economy, Keio University, Hiyoshi Hokuku-ku 
Yokohama, Japan, and Korea Institute for Advanced Study, 
207-43 Cheonryangni-2dong, Dongdaemun-gu, Seoul 130-722, Korea
}
\email{oguiso@hc.cc.keio.ac.jp}

\subjclass[2000]{14J50, 14J28}
\begin{abstract} 
Mordell-Weil groups of different abelian fibrations of a hyperk\"ahler manifold 
may have non-trivial relation even among elements of infinite order, 
but have essentially no relation, as its birational 
transformation. Precise definition of the terms 
"essentially no relation" will be given in Introduction. 
\end{abstract}

\maketitle


\setcounter{section}{0}
\section{Introduction - Background and the statement of main result} 

Unless stated otherwise, varieties are assumed to be normal, projective and defined over $\mathbf C$. Main results are Theorems (1.5) and (1.6).

Let $M$ be a hyperk\"ahler manifold, i.e. a simply-connected projective manifold admitting an everywhere non-degenerate global 
holomorphic $2$-form $\sigma_{M}$ s.t. $H^{0}(M, \Omega_{M}^{2}) = 
\mathbf C \sigma_{M}$ ([GHJ, Part III]). We denote by ${\rm Bir}\, M$ the 
group of birational transformations of $M$. A morphism 
$\varphi : M \longrightarrow B$ onto a variety $B$ is called an 
{\it abelian fibration} 
if its generic fiber $M_{\eta} (\ni O)$ is a positive dimensional abelian variety 
defined over $\mathbf C(B)$. The group $M_{\eta}(\mathbf C(B))$, 
called the Mordell-Weil group of $\varphi$ and denoted from now by ${\rm MW}(\varphi)$, 
can be naturally regarded as an abelian subgroup of ${\rm Bir}\, M$. In [Og3](see also [Og4]), 
we have shown the following:

\begin{theorem} \label{theorem:tits}

(1) Let $G < {\rm Bir}\, M$. Then either $G > \mathbf Z * \mathbf Z$ or $G$ is 
an almost abelian group of rank at most ${\rm Max}\, (1, \rho(M) -2)$, the latter 
of which is finitely generated. 

(2) Assume that $M$ admits $2$ different abelian fibrations $\varphi_{i} : M \longrightarrow B_{i}$ s.t. ${\rm rank}\, {\rm MW}(\varphi_{i}) > 0$ ($i =1$, $2$). Let $f_{i} \in {\rm MW}(\varphi_{i})$ s.t. ${\rm ord}\, f_{i} = \infty$. Then the subgroup $\langle f_{1}, f_{2} \rangle$ of ${\rm Bir}\, M$ contains $\mathbf Z * \mathbf Z$ as its subgroup. 
\end{theorem}  

Unfortunately, the proof in [Og3], which is based on Tits' alternative 
[Ti] together with properties of Salem polynomial (cf. [Mc]), tells us no explicit 
way to 
find $\mathbf Z * \mathbf Z$ in $\langle f_{1}, f_{2} \rangle$. Then F. Campana asked the following:

\begin{question} \label{question:campana} 

Under the same condition as in Theorem (1.1)(2), 
$$\langle f_{1}, f_{2} \rangle = \langle f_{1} \rangle * \langle f_{2} \rangle 
\simeq \mathbf Z * \mathbf Z$$
for {\it suitably chosen} $f_{i} \in {\rm MW}(\varphi_{i})$ ($i = 1$, $2$)? 
\end{question}

In the view of an observation of Cantat [Ca], one might ask  
an even stronger: 

\begin{question} \label{question:strong} 
Under the same condition as in Theorem (1.1)(2), 
$$\langle f_{1}, f_{2} \rangle = \langle f_{1} \rangle * \langle f_{2} \rangle 
\simeq \mathbf Z * \mathbf Z$$
for {\it any} $f_{i} \in {\rm MW}(\varphi_{i})$ s.t. ${\rm ord}\, (f_{i}) = \infty$ ($i = 1$, $2$) at least when ${\rm dim}\, M = 2$? 
\end{question}

The aim of this note is to give an affirmative answer to the first question in a slightly more general form and a negative answer to the second question. We note that there are infinitely many K3 surfaces and hyperk\"ahler manifolds admitting $2$ or more different abelian fibrations of positive Mordell-Weil 
rank (see eg. [Og3]). 

\begin{definition} \label{definition:norelation} 

Let $G$ be a group and let $G_{i}$ ($i =1, 2, \cdots, s$) be subgroups 
of $G$. We say that the groups $G_{i}$ ($i =1, 2, \cdots, s$) have 
{\it essentially no relation} if there are finite index 
subgroups $H_{i} < G_{i}$ s.t. 
$$\langle H_{1}, H_{2}, \cdots, H_{s} \rangle = 
H_{1} * H_{2} * \cdots * H_{s}\,\, .$$ 
\end{definition}

Our main result is the following:

\begin{theorem} \label{theorem:main} 

(1) Let $\varphi_{i} : M \longrightarrow B_{i}$ ($i = 1, 2, \cdots, s$) 
be mutually different abelian fibrations on a hyperk\"ahler manifold 
$M$ s.t. ${\rm rank}\, MW(\varphi_{i}) =: r_{i} > 0$. Then, 
in ${\rm Bir}\, M$, the Mordell-Weil groups 
${\rm MW}(\varphi_{i})$ ($i =1, 2, \cdots, s$) have essentially no relation, 
i.e. there are finite index subgroups $H_{i} < {\rm MW}(\varphi_{i})$ s.t. 
$$\langle H_{1}, H_{2}, \cdots, H_{s} \rangle = 
H_{1} * H_{2} * \cdots * H_{s} \simeq \mathbf Z^{r_{1}} * \mathbf Z^{r_{2}} * 
\cdots * \mathbf Z^{r_{s}}\,\, .$$ 
In particular, for given $g_{i} \in {\rm MW}(\varphi_{i})$ with 
${\rm ord}\, g_{i} = \infty$, there are positive integers $m_{i}$ 
($i=1, 2, \cdots, s$) s.t. 
$$\langle g_{1}^{m_{1}}, g_{2}^{m_{2}}, \cdots, g_{s}^{m_{s}} \rangle 
= \langle g_{1}^{m_{1}} \rangle * \langle g_{2}^{m_{2}} \rangle * 
\cdots * \langle g_{s}^{m_{s}} \rangle \simeq \mathbf Z * \mathbf Z * \cdots * 
\mathbf Z\,\,(s{\rm -factors})\,\, .$$

(2) There are a K3 surface $S$ admitting $2$ different Jacobian fibrations 
$\varphi_{i} : S 
\longrightarrow \mathbf P^{1}$ and elements $f_{i} \in {\rm MW}(\varphi_{i})$ 
with ${\rm ord}\, f_{i} = \infty$ ($i =1, 2$) s.t. in ${\rm Aut}\, S$, 
$$\langle f_{1}\, ,\, f_{2}\rangle \simeq \mathbf Z/2 * \mathbf Z/3\,\, , 
{\rm whence}\,\, ,\,\, \langle f_{1}\, ,\, f_{2}\rangle \not= \langle f_{1} \rangle * \langle f_{2} 
\rangle (\simeq \mathbf Z * \mathbf Z)\,\, ,$$ 
but, for each integer $n \ge 2$, 
$$\langle f_{1}^{n}\, ,\, f_{2}^{n}\rangle = \langle f_{1}^{n} \rangle * \langle f_{2}^{n} 
\rangle \simeq \mathbf Z * \mathbf Z\,\, .$$ 
\end{theorem}

Recall that a K3 surface is a $2$-dimensional hyperk\"ahler manifold, and an abelian fibration on a 
surface is usually called a Jacobian fibration. We also note that ${\rm Aut}\, S = {\rm Bir}\, S$ for a K3 surface (or more generally 
for a minimal surface). We shall prove Theorem (1.5)(1) in Section 2 and Theorem (1.5)(2) in Section 3. 
Our proof here is much inspired by an elementary hyperbolic geometry explained in the book [Th] of Thurston and an account [Ha] of de la Harpe, whose applicability is suggested by S. Cantat in connection with [Og3]. 

As a more concrete example, we shall also show, in Section 4, the following:

\begin{theorem} \label{theorem:k3} There is a K3 surface $S$ s.t. for each 
positive integer $s$, 
$${\rm Aut}\, S > \mathbf Z^{18} * \mathbf Z^{18} * \cdots * \mathbf Z^{18}\,\, (s{\rm -factors})\,\, .$$
\end{theorem}

Note that if $\mathbf Z^{r} < {\rm Aut}\, S$ for some K3 surface $S$, then 
$r \le 18$ by Theorem (1.1)(2).

{\it Acknowledgement.} I would like to express my deep thank 
to Professors F. Campana, S. Cantat, J.M. Hwang, M. Paun, and J. Winkelman 
for their 
valuable discussions, mostly in Johns Hopkins University (Recent Developments in Higher Dimensional Algebraic Geometry, March 2006), Luminy (Dynamics and Complex Geometry, June 2006) and Oberwolfach (Komplexe Analysis, August 2006).

\section{Proof of Theorem (1.5)(1)}

{\bf 1.} The next Theorem and its proof are taken from Tits' article [Ti]. This simple theorem 
is very useful in our proof of Theorem (1.5).

\begin{theorem} \label{theorem:pingpong} Let $s$ be an integer s.t. $s \ge 2$ 
and let $G$ be a group which acts faithfully on a non-empty set $S$. Let $G_{i}$ ($i = 1, 2, \cdots, s$) be subgroups of $G$, let $S_{i}$ 
($i = 1, 2, \cdots, s$) be subsets of $S$ s.t. $S \setminus \cup_{i=1}^{s} S_{i} \not= \emptyset$, and let $p$ be an element of $S \setminus \cup_{i=1}^{s} S_{i}$. Put $G_{i}^{0} := G_{i} \setminus \{id \}$ and $\tilde{S}_{i} := S_{i} \cup \{p\}$. Suppose that $G_{i}^{0}(\tilde{S}_{j}) \subset S_{i}$ for all pairs $i \not= j$. Then 
$$\langle G_{1}, G_{2}, \cdots , G_{s} \rangle = G_{1} * G_{2} * \cdots * G_{s}\,\, .$$
\end{theorem} 

\begin{proof} Consider $g_{1}g_{2} \cdots g_{m}$, where 
$g_{k} \in G_{i_{k}}^{0}$ and $i_{1} \not= i_{2} \not= \cdots \not= i_{m}$. 
Using $s \ge 2$, one can choose an integer $i$ s.t. $i \not= i_{m}$ and 
$1 \le i \le s$. Then $G_{i_{m}}^{0}(\tilde{S}_{i}) \subset S_{i_{m}}$, 
$G_{i_{m-1}}^{0}G_{i_{m}}^{0}(\tilde{S}_{i}) \subset S_{i_{m-1}}$, 
and finally $G_{i_{1}}^{0}G_{i_{2}}^{0}\cdots G_{i_{m-1}}^{0}G_{i_{m}}^{0}(\tilde{S}_{i}) \subset S_{i_{1}}$. Thus $g_{1}g_{2} \cdots g_{m}(p) 
\subset S_{i_{1}}$. 
Hence $g_{1}g_{2} \cdots g_{m}(p) \not= p$ by $p \not\in S_{i_{1}}$, 
and therefore $g_{1}g_{2} \cdots g_{m} \not= id$. 
\end{proof} 

{\bf 2.} Let $M$ be a hyperk\"ahler manifold. As it is mentioned in 
Introduction, $M$ is assumed to be projective. In this paragraph, we 
recall some basic facts about $M$ and fix some notations.

{\bf (2.1)} The N\'eron-Severi group $NS(M)$ is a non-degenerate 
lattice of singature $(1, \rho(M) -1)$ w.r.t. the Beauville form 
$(*,**)$. (See eg. [GHJ, Part III]; the Beauville form here is called 
the Beauville-Bogomolov form there). Here $\rho(M) := {\rm rank}\, NS(M)$ 
is the Picard number of $M$. Note that $c_{1} : 
{\rm Pic}\, M \longrightarrow NS(M)$ is an isomorphism for $M$. 

Let ${\rm O}(NS(M))$ be the group of isometries of the lattice $NS(M)$. 
Then, there is a natural group homomorphism 
$$r_{NS} : {\rm Bir}\, M \longrightarrow {\rm O}(NS(M))\,\, ;\,\, 
g \mapsto (g^{-1})^{*} \vert NS(M)\,\,$$
(see eg. [GHJ, Part III]), and $\vert {\rm Ker}\, r_{NS} \vert < \infty$ 
(see eg. [Hu], [Og2]). For the last statement, 
projectivity of $M$ is essential ([Og2]).  

{\bf (2.2)} In what follows, we denote the scalar extension 
$NS(M) \otimes K$ ($K$ is a field) simply by $NS(M)_{K}$. Whenever we discuss 
topology, we regard $NS(M)_{\mathbf R} \simeq \mathbf R^{\rho(M)}$ as the Euclidean space by a suitable Euclidean norm $\vert\vert * \vert\vert$. 

Let $\mathcal C$ be the {\it positive cone} of $NS(M)$, that is, 
the connected component of $\{ x \in NS(M)_{\mathbf R} \vert (x, x) > 0\}$ 
containing ample classes. We denote by $\overline{\mathcal C}$ the closure 
of $\mathcal C$ in $NS(M)_{\mathbf R}$, and by 
$\partial \overline{\mathcal C}$, 
the boundary of $\mathcal C$. We define
$$S := \{x \in \partial \mathcal C\, \vert\,\,  \vert\vert x \vert\vert = 1\}\,\, .$$
This space $S$ can be identified with the set of the rays 
$\{\mathbf R_{+}x\, \vert\, x \in \partial \mathcal C \setminus \{0\}\}$, by the natural bijection $\mathbf R_{+}x \leftrightarrow x/\vert\vert x\vert\vert$, and is also 
homeomorphic to the sphere $S^{\rho(M) -2}$. 

{\bf (2.3)} Let $\varphi : M \longrightarrow B$ be an abelian fibration on 
$M$. Then, by a fundamental result due to Matsushita [Ma],  $B$ is 
$\mathbf Q$-factorial and of Picard number $1$. 
Let $v$ be the (uniquely determined) primitive element in $\mathbf R_{+} \varphi^{*}h_{B} \cap NS(M)$. Here $h_{B}$ is an ample class of ${\rm Pic}\, B$. 
We call this $v$ the {\it primitive vector corresponding to} $\varphi$. Note 
that $v \in (\partial \overline{\mathcal C} \cap NS(M)) \setminus \{0\}$. We denote $v/\vert\vert v \vert\vert \in S$ by $\overline{v}$ and call $\overline{v}$ the {\it unit vector corresponding to} $\varphi$.

Let $\tilde{G}$ be one of the maximal free abelian subgroups of 
$${\rm Im}(r_{NS} \vert {\rm MW}(\varphi) : {\rm MW}(\varphi) \longrightarrow 
{\rm O}(NS(M)))\,\, .$$

Then $\tilde{G} \simeq \mathbf Z^{r}$ ($r = {\rm rank}\, {\rm MW}(\varphi)$)  
by the paragraph (2.1). Note that the primitive vector $v$ is stable under $\tilde{G}$, i.e. 
$\tilde{G}(v) = v$. This is because $\rho(B) = 1$. 

Given another abelian fibration $\varphi' : M \longrightarrow B'$, we write 
the primitive (resp. unit) element of $\varphi'$ by $v'$ 
(resp. by $\overline{v'}$). Again by $\rho(B) = \rho(B') = 1$, we have $\varphi = \varphi'$ iff $v = v'$, or equivalently, iff $\overline{v} = \overline{v'}$. 

{\bf (2.4)} Let $\varphi : M \longrightarrow B$, $v$, $\overline{v}$, 
$\tilde{G}$ be the same as in (2.3). Following [Th, Problem 2.5.24], 
we consider the spaces:
$$\tilde{L}^{\mathbf Q} := \{x \in NS(M)_{\mathbf Q} \vert (x, v) = 1\}\,\, ,\,\, 
\tilde{L} := \{x \in NS(M)_{\mathbf R} \vert (x, v) = 1\}\,\, ;$$
$$L^{\mathbf Q} := \tilde{L}^{\mathbf Q}/\mathbf Q v\,\, ,\,\, 
L := \tilde{L}/ \mathbf R v\,\, ;$$
$$P := \{x \in \partial \mathcal C\, \vert\, (x, v) = 1\} \subset \tilde{L}\,\, .$$
Here, since $v \in NS(M)$ and $(v, v) =0$, 
the mapping $av : x \mapsto x + av$, $a \in \mathbf Q$ (resp. 
$a \in \mathbf R$) 
defines a faithful action of $\mathbf Q v$ 
(resp. of $\mathbf R v$) on $\tilde{L}^{\mathbf Q}$ (resp. on $\tilde{L}$). 
We note that these spaces depend on (the choice of) an abelian fibration 
$\varphi : M \longrightarrow B$. 

\begin{remark} \label{remark:thurston}Though the rational spaces 
$\tilde{L}^{\mathbf Q}$ and $L^{\mathbf Q}$ are not considered in [Th, Problem 2.5.14], 
it is important to consider these spaces for our proof of Theorem (1.5). 
\end{remark}

By definition, the space $L^{\mathbf Q}$ (resp. $L$) is a $(\rho(M) -2)$-dimensional rational 
(resp. real) affine space. Or more explicitly, 
by taking a point in $u \in L^{\mathbf Q}$ as $0$ 
and by taking a rational basis of the tangent space at $u$, we can identify 
$L^{\mathbf Q}$ (resp. $L$) with the vector space $\mathbf Q^{\rho(M) -2}$ 
(resp. $\mathbf R^{\rho(M) -2}$). Under these identifications, we have also 
$L = (L^{\mathbf Q})_{\mathbf R}$. We denote by $p$ the natural 
quotient maps :   
$$p : \tilde{L}^{\mathbf Q} \longrightarrow L^{\mathbf Q}\,\, ,\,\, 
p : \tilde{L} \longrightarrow L\,\, .$$ 

Finally, we define the map $\pi$ by 
$$\pi : P \longrightarrow S\,\, ;\,\, x \mapsto x/\vert\vert x \vert\vert
\,\, .$$

Since $\tilde{G}v = v$ and $(v, v) = 0$, the group $\tilde{G}$ acts on $\tilde{L}^{\mathbf Q}$, $\tilde{L}$, $L^{\mathbf Q}$, $L$, $P$ and $S$. Here 
the action of $\tilde{G}$ on $S$ is defined by 
$x \mapsto g(x)/\vert\vert g(x) \vert\vert$. By definiton, 
the actions of $\tilde{G}$ on these spaces are equivariant under $p$ and $\pi$.

{\bf 3.} In this paragraph, again closely following [Th, Problem 2.5.14], 
we shall prove three propositions, which are crucial in our proof of Theorem 
(1.5)(1). Let $\varphi : M \longrightarrow B$ be the same as in the paragraph 
2. We use the same notation as in the paragraphs (2.3) and (2.4). 

\begin{proposition} \label{proposition:compactif} 

The map $p \vert P : P \longrightarrow L$ is a homeomorphism and the map  
$$\iota := \pi \circ (p \vert P)^{-1} : L \longrightarrow S$$ 
is a homeomorphism onto $S \setminus \{\overline{v}\}$, i.e. 
$S$ is a one point compactification of the real affine space $L$ through 
$\iota$. Moreover, the actions of $\tilde{G}$ on $L$ and $S$ are equivariant 
under $\iota$.  
\end{proposition} 

\begin{proof} The last assertion follows from the last statement in the paragraph (2.4). Let us show the first two assertions. For a given $x \in \tilde{L}$, 
there is a unique $\alpha$ s.t. $(x + \alpha v, x + \alpha v) = 0$, namely $\alpha = -(x, x)/2$. Since $(x + \alpha v, v) = 1 > 0$, we have $x + \alpha v \in 
\partial \mathcal C$, and therefore $x + \alpha v \in P$. By the uniqueness of $\alpha$, 
the map $x + \mathbf R v \mapsto x + \alpha v$ is well-defined 
and gives the inverse $q : L \longrightarrow P$ of $p \vert P$. 

By $(v, v) = 0$, we have $\pi(P) \subset S \setminus \{\overline{v}\}$. 
Let $x \in \partial \mathcal C \setminus \{0\}$. By the Schwartz inequality, 
we have $(x, v) > 0$ unless $x \in \mathbf Rv$. Thus, for $x \in S \setminus \{\overline{v}\}$, we have $x/(x,v) \in P$ and $\pi(x/(x,v)) = x$. Thus, 
$x \mapsto x/(x, v)$ gives the inverse of $\pi : P \longrightarrow S \setminus \{\overline{v}\}$. Hence $\iota = \pi \circ q$ satisfies the requirement. 
\end{proof}

\begin{proposition} \label{proposition:parallel} 
There is a finite index subgroup $G$ of $\tilde{G}$ s.t. 
$G$ acts on the real affine space $L$ as parallel transformations. 
Moreover, the action of $G$ on $L$ is faithful and discrete. 
\end{proposition} 

\begin{remark} For both statements, it seems quite essential that the 
action $\tilde{G}$ is defined over the {\it lattice} $NS(M)$, or in other 
words, over $\mathbf Z$. 
Being parallel transformations is false 
even for a cyclic orthogonal subgroup of $NS(M)_{\mathbf Q}$ fixing 
some non-zero rational vector $v$ with $(v,v) = 0$, as there is an 
infinite-order $3 \times 3$ rotation matrix of rational entries. 
An account of [Ha, Page 135] seems to miss this point. See remark (2.7) 
for discreteness. 
\end{remark} 

\begin{proof} Since $v$ is a primitive vector with $(v, v) = 0$, one can choose an integral basis of $v^{\perp}_{NS(M)}$ as $\langle v\, ,\, w_{1}\, ,\, \cdots \, ,\, w_{n}\,\rangle$. 

Choose $u \in NS(M)_{\mathbf Q}$ s.t. $(v, u) = 1$. One can choose such $u$, 
because $NS(M)$ is non-degenerate. Then 
$$\langle\, v\,\, ,\,\, w_{1}\,\, ,\,\, \cdots \,\, ,\,\, w_{n}\,\, ,\,\, u\, 
\rangle$$
forms $\mathbf Q$-basis of $NS(M)_{\mathbf Q}$. Note that, 
in this notation, $\rho(M) = n +2$. 

Let $g \in \tilde{G}$. Since $g(v) = v$ and $g(v^{\perp}_{NS(M)}) = 
v^{\perp}_{NS(M)}$, the matrix representation $M(g)$ of $g$ w.r.t. 
the basis above is of the form:

$$M(g) :=\, \left( \begin{array}{ccc} 1& ^{t}\mathbf a(g) & c(g)\\
\mathbf o &  A(g) & \mathbf b(g)\\
0 & ^{t}\mathbf o & d(g) \end{array}\right)\, .$$ 
Here $A(g)$ is the matrix representation of the action of $g$ on 
the {\it lattice} $N := v^{\perp}_{NS(M)}/\mathbf Z v$ w.r.t. its 
{\it integral} basis $\langle [w_{i}] \rangle_{i=1}^{n}$. Thus $A(g)$ is 
an {\it integral} matrix, while $\mathbf b(g) \in \mathbf Q^{n}$ is, in general, a rational vector.

Since $N$ is of negative definite, its orthogonal 
transformation $A(g)$ is diagonalizable, and the eigenvalues of 
$A(g)$ are of absolute value $1$. On the other hand, 
since $A(g)$ is an integral matrix, its eigenvalues are algebraic 
integers. Thus, by Kronecker's theorem, the eigenvalues of $A(g)$ 
are roots of $1$. Since $A(g)$ is diagonalizable, $A(g)$ is then 
of finite order. Let us denote by $m(g)$ the order of $A(g)$.

Let $\langle {\tilde g}_{i} \rangle_{i=1}^{r}$ be a generator of the free 
abelian group 
$\tilde{G}$. Put $g_{i} := {\tilde g}_{i}^{2m(\tilde g_{i})}$ and define 
the subgroup $G$ of $\tilde{G}$ by
$$G := \langle g_{1}\, ,\, g_{2}\, ,\, \cdots \, ,\, g_{r} \rangle\,\, .$$

Since $\tilde{G}$ is a finitely generated abelian group, 
$G$ is a finite index subgroup of $\tilde{G}$. Moreover, for each 
$g \in G$, 
we have $A(g) = I_{n}$ (the identity matrix) and $d(g) = 1$, i.e. 
$$M(g) :=\, \left( \begin{array}{ccc} 1& ^{t}\mathbf a(g) & c(g)\\
\mathbf o &  I_{n} & \mathbf b(g)\\
0 & ^{t}\mathbf o & 1 \end{array}\right)\,\, .$$ 

For $d(g) = 1$, observe that ${\rm det}\, M(\tilde{g}_{i}) = \pm 1$, whence 
${\rm det}\, M(\tilde{g}_{i}^{2m(\tilde{g}_{i})}) = 1$, and therefore 
${\rm det}\, M(g) = 1$ for $g \in G$. 

Let $Q \in L$. Then $Q$ is uniquely expressed in the form: 
$$Q = [u] + \sum_{i=1}^{n} x_{i}(Q)[w_{i}]\,\, ,$$
and the vector valued function $\mathbf x := (x_{i})_{i=1}^{n}$ gives an affine coordinate of both $L$ and $L^{\mathbf Q}$. As it is remarked in the 
paragraph (2.4), under this coordinate, one can identify 
$$L^{\mathbf Q} = \mathbf Q^{n}\, ,\, L = (\mathbf Q^{n}) \otimes \mathbf R = \mathbf R^{n}\,\, .$$ 
By the shape of $M(g)$, the action of $g \in G$ on $L = \mathbf R^{n}$ is 
expreessed by the parallel transformation by $\mathbf b(g) \in \mathbf Q^{n}$. 
This shows the first statement. 

For the discreteness and faithfulness of the action of $G$ on 
$L = \mathbf R^{n}$, it suffices to show Lemma (2.6) below. Indeed, 
$\langle \rho(g_{i})\rangle_{i=1}^{r}$ then forms a part of 
real basis of $L = (\mathbf Q^{n}) \otimes \mathbf R$.

\begin{lemma} \label{lemma:inj} 
The map 
$$\rho : G \simeq \mathbf Z^{r} \longrightarrow \mathbf Q^{n}\,\, ;\,\, 
g \mapsto \mathbf b(g)$$
is an injective group homomorphism.  
\end{lemma} 

\begin{remark} Note that the map $(a, b) \mapsto a + b\sqrt{2}$ 
defines an injective group homomorphism from $\mathbf Z^{2}$ to 
the group $\mathbf R$ of the parallel transformations of $\mathbf R$, 
but the image $\{a + b\sqrt{2} \vert a, b \in \mathbf Z\}$ is not discrete 
in $\mathbf R$. So, the fact that $G$ acts on the underlying rational space 
$L^{\mathbf Q}$ seems crucial for the discreteness.
\end{remark}

\begin{proof} Assume that $\mathbf b(g) = \mathbf o$. Then, by the shape 
of $M(g)$, we have $g(v) = v$, $g(w_{i}) = w_{i} + a_{i}(g)$, $g(u) = u + c(g)v$. From 
$$(u, u) = (g(u), g(u)) = (u, u) + 2c(g)(u, v) = (u, u) + 2c(g)\,\, ,$$
we have $c(g) = 1$. Then, from 
$$(u, w_{i}) = (g(u), g(w_{i})) = (u, w_{i}) + a_{i}(g)(u, v) = (u, w_{i}) + a_{i}(g)\,\, ,$$
we have $a_{i}(g) = 0$. Hence $M(g) = I_{n+2}$, i.e. $g = id$. 

\end{proof}
This completes the proof of Proposition (2.4).
\end{proof}

\begin{proposition} \label{proposition:compact} 
Let $U \subset S$ be a compact neighborhood of the unit 
vector $\overline{v}$ and let $V$ be a non-empty compact subset of 
$S \setminus \{\overline{v}\}$. (See paragraphs (2.3) and (2.4) 
for the definition 
of $S$ and $\overline{v}$.) 
Then, there is a finite index subgroup $H$ of $\tilde{G}$ s.t. 
$(H \setminus \{id\})(V) \subset U$. 
\end{proposition} 

\begin{proof} Identifying $L$ with $S \setminus \{\overline{v}\}$ by $\iota$ (see Proposition (2.3)), we may assume that 
$$V \subset L\,\, {\rm and}\,\, B:= \{x \in L \vert\, \vert\vert x \vert\vert > r\} \subset U \setminus \{\overline{v}\} \subset L$$ 
for some $r > 0$. Here 
$\vert\vert * \vert\vert$ is some Euclidean norm of $L$ (w.r.t. some fixed origin). 

Let $G \simeq \mathbf Z^{r}$ be a finite index subgroup of $\tilde{G}$ found in Proposition (2.4). 
Then, $G$ acts on $L$ by parallel transformations, say by 
$\{t(g) \in \mathbf R^{n} \vert g \in G\}$. 

Since the action of $G$ is discrete by Proposition (2.4), we have

$${\rm inf}\, \{\vert\vert t(g) \vert\vert\,\, \vert\,\, g \in G \setminus \{id\}\} > 0\,\, .$$

Thus, there is some positive integer $c$ s.t. $g^{c}(V) \subset B$ for all $g \in G \setminus \{id\}$. Set $H := \{g^{c} \vert g \in G\}$. Since $G$ is a finitely generated abelian group, this set $H$ is a finite index subgroup of $G$, 
and therefore, is a finite index subgroup of $\tilde{G}$. This $H$ satisfies the requirement. 
\end{proof}

{\bf 4.} Now we shall complete the proof of Theorem (1.5)(1). 

Let $\varphi_{i} : M \longrightarrow B_{i}$ ($i = 1, \cdots , s$) be mutually 
different abelian fibrations. For each $i$, we denote by 
$\overline{v}_{i} \in S$ the unit element 
corresponding to $\varphi_{i}$. Note that 
$\overline{v}_{i} \not= \overline{v}_{j}$ if $i \not= j$. 

Put $r_{i} := {\rm rank}\, {\rm MW}\, (\varphi_{i})$ and choose 
one of the maximal free abelian subgroups $\tilde{G}_{i} \simeq \mathbf Z^{r_{i}}$ of  
${\rm Im}\, (r_{NS} : {\rm MW}(\varphi_{i}) 
\longrightarrow {\rm O}(NS(M)))$ for each $i$. 

Let $U_{i} \subset S$ be a compact neighborhood of $\overline{v}_{i}$. 
We can take $U_{i} \cap U_{j} = \emptyset$ if $i \not= j$ 
and $S \setminus \cup_{i=1}^{s} U_{i} \not= \emptyset$. 
Choose $p \in S \setminus \cup_{i=1}^{s} U_{i}$ 
and put $\overline{U}_{i} := U_{i} \cup \{p\}$. Then, by Proposition (2.8), 
there are finite index subgroups $H_{i}'$ of $G_{i}$ s.t. 
$$(H_{i}' \setminus \{id\})(\overline{U}_{j}) \subset U_{i}$$ 
for all pairs $i \not= j$ in $\{1, 2, \cdots, s\}$. Then, by Theorem (2.1), 
we obtain 
$$\langle H_{1}', H_{2}', \cdots , H_{s}' \rangle = H_{1}' * H_{2}' * \cdots 
* H_{s}' 
\simeq \mathbf Z^{r_{1}} * \mathbf Z^{r_{2}} * \cdots * \mathbf Z^{r_{s}}\,\, 
.$$ 
One can choose free abelian subgroups $H_{i}$ of ${\rm MW}(\varphi_{i})$ s.t. 
$r_{NS}(H_{i}) = H_{i}'$. These $H_{i}$ satisfy the requirement.

\section{Proof of Theorem (1.5)(2)}

Let $E$ be an elliptic curve. Let us consider the product abelian surface 
$A := E \times E$, and its associated Kummer K3 surface $S := {\rm Km}\, A$. 
By definition, $S$ is the minimal resolution of the quotient surface $A/\iota$ ($\iota$ is defined below). 

Let $\tilde{p}_{i} : A \longrightarrow E$ be the projection to the 
$i$-th factor, 
and let $p_{i} : S \longrightarrow \mathbf P^{1}$ be the Jacobian fibration 
on $S$, induced by $\tilde{p}_{i}$. 

Let 

$$\tilde{f}_{2} := \left( \begin{array}{cc} 1& 1\\
0& 1 \end{array}\right)\,\, ,\,\, 
\tilde{f}_{1} := \left( \begin{array}{cc} 1& 0\\
1& 1 \end{array}\right)\,\, ,\,\, \iota :=  \left( \begin{array}{cc} -1& 0\\
0& -1 \end{array}\right)$$
and put $\tilde{G} := \langle \tilde{f}_{1}, \tilde{f}_{2} \rangle$ in ${\rm SL}\, (2, \mathbf Z)$. Then $\tilde{G} = {\rm SL}(2, \mathbf Z)$ 
(see eg. [Kn, Page 80 Exercise (3.2)]), 
and ${\rm SL}(2, \mathbf Z)$ acts faithfully on $A$ by 
$$\tilde{f}_{2}(x, y) = (x+y, y)\,\, ,\,\, 
\tilde{f}_{1}(x, y) = (x, x+y)\,\, .$$ 
By the shape of $\tilde{f}_{i}$, we see that 
 $\tilde{f}_{i} \in {\rm MW}\, (\tilde{p}_{i})$ and ${\rm ord}\, \tilde{f}_{i} 
= \infty$. Since $\tilde{f}_{i}\iota = 
\iota \tilde{f}_{i}$, we also see that $\tilde{f}_{i}$ descends to $f_{i} \in {\rm MW}(p_{i})$ and ${\rm ord}\, f_{i} = \infty$. 

Set $G := \langle f_{1}, f_{2} \rangle$ in ${\rm Aut}\, S$. Then, 
$G \simeq {\rm SL}(2, \mathbf Z)/\langle \iota \rangle = {\rm PSL}(2, \mathbf Z)$. 

We shall show that $(S, f_{1}, f_{2})$ satisfies the requirement. 

It is well-known that ${\rm PSL}(2, \mathbf Z) \simeq 
\mathbf Z/2 * \mathbf Z/3$ (see eg. [Kn, Page 147]). In particular, 
$\langle f_{1}, f_{2} \rangle \not= \langle f_{1} \rangle * \langle 
f_{2} \rangle$; if otherwise, $\mathbf Z/2 * \mathbf Z/3 
\simeq \mathbf Z * \mathbf Z$, a contradiction. 

It remains to show that $\langle f_{1}^{n}, f_{2}^{n} \rangle = \langle 
f_{1}^{n} \rangle * \langle f_{2}^{n} \rangle$ for each $n \ge 2$. Our proof 
of this fact closely follows [Ha, Example 1]. Note that 
${\rm PSL}(2, \mathbf Z) < {\rm PGL}(2, \mathbf Z)$. 

In ${\rm PGL}(2, \mathbf Z)$, we put 
$$g_{2} := f_{2}^{n} = \left( \begin{array}{cc} 1& n\\
0& 1 \end{array}\right)\,\, ,\,\, 
g_{1} := f_{1}^{n} := \left( \begin{array}{cc} 1& 0\\
n& 1 \end{array}\right)\,\, ,\,\, j :=  \left( \begin{array}{cc} 0& 1\\
1& 0 \end{array}\right)\,\, .$$
Then $j^{2} = id$ and $g_{1} = jg_{2}j$ in ${\rm PGL}(2, \mathbf Z)$. So, 
if there would be a non-trivial relation among $g_{1}$ and $g_{2}$, say, $h(g_{1}, g_{2}) = id$, then substituting $g_{1}^{l} = jg_{2}^{l}j$ into the relation 
$h(g_{1}, g_{2}) = id$, we would have a non-trivial relation among $g_{2}$ and $j$. Thus, we suffice to show 
that $\langle g_{2}, j \rangle = \langle g_{2} \rangle * \langle j \rangle$ 
in ${\rm PGL}(2, \mathbf Z)$. 

To prove this claim, we consider the natural fractional linear action of 
${\rm PGL}(2, \mathbf Z)$ on $\mathbf P^{1} = \mathbf C \cup \{\infty\}$, 
and the following 
subsets and a point of $\mathbf P^{1}$:
$$U_{1} := \{z \in \mathbf C\, \vert\, \vert z \vert < 1\}\,\, ,\,\, 
U_{2} := \{z \in \mathbf C\, \vert\, \vert {\rm Re}\, z \vert > 1\} 
\cup \{\infty\}\,\, ,\,\, P := 2\sqrt{-1}\ \not\in U_{1} \cup U_{2}\,\, .$$

Then $j(U_{2} \cup \{P\}) \subset U_{1}$ and 
$g_{2}^{k}(U_{1} \cup \{P\}) \subset U_{2}$ for each 
$k \not= 0$ (by $n \ge 2$). Thus, by Theorem (2.1), we have 
$\langle g_{2}\, ,\, j \rangle = \langle g_{2} \rangle * \langle j \rangle$, 
and we are done. 

\section{Proof of Theorem (1.6)} 

Let $\varphi : S \longrightarrow \mathbf P^{1}$ be a Jacobian K3 surface s.t. 
${\rm rank}\, {\rm MW}(\varphi) = 18$ and $\rho(S) = 20$. 
Such a Jacobian K3 surface exists by 
[Co] and [Ni]. It is also known that, given a Jacobian K3 surface $M$, 
one can find such 
a Jacobian K3 surface arbitrarily close to $M$ in the period domain ([Og1]). 

We shall show that this $S$ satisfies the requirement. 

Let $e \in {\rm Pic}\, S 
\simeq NS(S)$ be the fiber class of $\varphi$. Since $\varphi$ has at least 
three singular 
fibers (see eg. [Ca]), the subgroup 
$${\rm Aut}\, \varphi := \{f \in {\rm Aut}\, S\, \vert\, f^{*}(e) = e\}$$
of ${\rm Aut}\, S$ is an almost abelian group (of rank $18$).

On the other hand, since $\rho(S) = 20$, we have $\mathbf Z * \mathbf Z \subset {\rm Aut}\, S$ by 
[Og3]. In particular, ${\rm Aut}\, S$ is not an almost abelian group. 

Hence the coset ${\rm Aut}\, S/{\rm Aut}\, \varphi$ is an infinite set. 
Thus, the ${\rm Aut}\, S$-orbit of $e$ is an infinite set as well. Therefore, 
$S$ admits infinitely many different Jacobian fibrations with Mordell-Weil rank $18$. 
Now the result follows from Theorem (1.5)(1).  


\end{document}